\documentclass[leqno,12pt]{amsart}
\usepackage{latexsym}
\usepackage{mathrsfs}
\usepackage{amsmath}
\usepackage{amssymb}
\usepackage[bookmarks]{hyperref}
\usepackage{pstricks,pst-plot}

\newtheorem{theorem}{Theorem}[section]
\newtheorem{lemma}[theorem]{Lemma}
\newtheorem{corollary}[theorem]{Corollary}
\newtheorem{proposition}[theorem]{Proposition}

\newtheorem{conjecture}[theorem]{Conjecture}

\theoremstyle{definition}

\newtheorem{definition}[theorem]{Definition}

\newcommand{\C}{\mathbb{C}}

\newcommand{\R}{\mathbb{R}}

\newcommand{\bthm}{\begin{theorem}}
\newcommand{\ethm}{\end{theorem}}
\newcommand{\blem}{\begin{lemma}}
\newcommand{\elem}{\end{lemma}}
\newcommand{\bcor}{\begin{corollary}}
\newcommand{\ecor}{\end{corollary}}
\newcommand{\bprop}{\begin{proposition}}
\newcommand{\eprop}{\end{proposition}}
\newcommand{\bdefn}{\begin{definition}}
\newcommand{\edefn}{\end{definition}}
\newcommand{\bpf}{\begin{proof}}
\newcommand{\epf}{\end{proof}}

\def\vep {\varepsilon}
\def \sm {\setminus}

\newcommand{\ra}{\rightarrow}

\def\row#1#2{#1_1,\ldots, #1_{#2}}

\def\rn{\R^n}
\def\u{\mathscr U}
\def\s{\sigma}
\def\t{\tau}
\def\th{\tilde h}
\def\tf{\tilde f}
\def\vv{\mathscr V}
\def\pv{\varphi_V}
\def\tg{\tilde g}
\def\tE{\widetilde E}
\def\tvep{\tilde \vep}
\def\talpha{\tilde\alpha}

\def\tgamma{\tilde\gamma}

\def\cpi{C^+_i}
\def\cmi{C^-_i}

\begin{document}
\title[Approximation by zero-free continuous maps]{Approximation by\\ zero-free continuous maps}

\author{Alexander J. Izzo}
\address{Department of Mathematics and Statistics, Bowling Green State University, Bowling Green, OH 43403}
\email{aizzo@bgsu.edu}

\begin{abstract}
We prove that if $E$ a subset of an $n$-dimensional manifold, then every continuous 
$\rn$-valued map on $E$ that is zero-free on the interior of $E$ can be approximated in the fine topology, and hence, in particular, in the uniform topology, by a continuous 
$\rn$-valued map that is zero-free on all of $E$.
\end{abstract}

\maketitle

%
%
%
%

\section{The Result}

The main purpose of this paper is to prove the following theorem.

\bthm\label{theorem}
Let $E$ be a subset of $\rn$, and let $f: E \ra \rn$ be a continuous map that is zero-free on the interior $E^\circ$ of $E$.  Then $f$ can be approximated uniformly by a
continuous $\rn$-valued map that is zero-free on all of $E$.
\ethm

Note that the subset $E$ of $\rn$ is completely arbitrary.

One would expect Theorem~\ref{theorem} to hold more generally for $E$ a subset of a manifold, and that is indeed the case.  As far as the author can see though, that more general result does not follow directly from the Euclidean case.  We will show, however, that the conclusion of Theorem~\ref{theorem} can be strengthened to yield a result that does transfer to the manifold setting.

Our most general result is as follows.

\bthm\label{theorem-generalized}
Let $E$ be a subset of an $n$-dimensional manifold $M$, and let $f: E \ra \rn$ be a continuous map that is zero-free on the interior $E^\circ$ of $E$ in $M$.  Then $f$
can be approximated in the fine topology by a 
continuous $\rn$-valued map that is zero-free on all of $E$.
\ethm

Here and throughout the paper, by an \emph{$n$-dimensional manifold} we mean a paracompact Hausdorff space every point of which has a neighborhood that is homeomorphic to an open subset of $\rn$.  The \emph{fine topology} on the set $C(E,\rn)$ of continuous maps from $E$ to $\rn$ is the topology determined by taking as a basis all sets of the form 
$$B(f,\vep)=\{g\in C(E,\rn): |f(x)-g(x)|<\vep(x) {\rm\ for\ all\ } x\in E\},$$
for $f\in C(E,\rn)$ and $\vep: E\ra \R_+=\{t\in \R: t>0\}$ a positive continuous function.
Thus explicitly, the conclusion of Theorem~1.2 is that for every positive continuous function $\vep: E\ra \R_+$ there is a map $g\in C(E,\rn)$ that is zero-free and satisfies $|f(x)-g(x)|<\vep(x)$ for every $x\in E$.

Note that Theorem~\ref{theorem} is an immediate consequence of Theorem~\ref{theorem-generalized}.

As the interested reader can verify, our proof of Theorem~\ref{theorem-generalized} carries over almost without change to show that in the theorem the word \lq\lq manifold\rq\rq\ can be replaced by \lq\lq manifold-with-boundary\rq\rq.

Our theorems are closely related to results in dimension theory.  It is standard that the topological dimension (i.e., Lebesgue covering dimension) $\dim E$ of a subset $E$ of an $n$-dimensional manifold $M$ satisfies $\dim E\leq n-1$ if and only if the interior of $E$ in $M$ is empty \cite[Theorem~III~1 and Corollary~1 of Theorems~IV~3]{Hurewicz-Wallman1948}.  It is also standard that every continuous $\rn$-valued map on a separable metrizable space $E$ is a uniform limit of zero-free $\rn$-valued maps if and only if $\dim E\leq n-1$ \cite[Theorems~VI~1 and~VI~2]{Hurewicz-Wallman1948}.  Together these two facts immediately yield the special case of Theorem~\ref{theorem} in which the interior of $E$ is empty.

The motivation for considering Theorem~\ref{theorem} comes from a conjecture in complex analysis.  Let $K$ be a compact set in the complex plane.  It is easily shown that if every continuous complex-valued function on $K$ that is holomorphic on $K^\circ$ is a uniform limit of polynomials in the complex coordinate function $z$, then the complement $\C\sm K$ of $K$ must be connected.  A famous theorem of S.~N.~Mergelyan \cite{Mer} (see also \cite[Theorem~20.5]{Rudin}) asserts that conversely, if the complement of $K$ is connected, then every continuous complex-valued function on $K$ that is holomorphic on $K^\circ$ is a uniform limit of polynomials in $z$.  A conjecture of Johan Andersson and Paul Gauthier asserts that when the function to be approximated is zero-free on $K^\circ$, then the approximating functions can be taken to be zero-free on all of $K$.

\begin{conjecture}[Andersson-Gauthier~\cite{AG}]
If $\mathbb C\setminus K$ is connected, then every continuous function $f:K\to \mathbb C$ that is holomorphic and zero-free on $K^\circ$ is a uniform limit of  polynomials that are zero-free on all  of $K$. 
\end{conjecture}

Theorem~\ref{theorem} show that there is no topological obstruction to this conjecture.

The special case of Theorem~\ref{theorem} in which the set $E$ is compact and $n=2$, which is the case relevant to the conjecture of Andersson and Gauthier, was proved by Christoper Bishop and Gauthier (private communication).  Although the proof we give for our theorems is different from their proof, both the statement of Theorem~\ref{theorem} and the proof we give of the more general Theorem~\ref{theorem-generalized} were inspired by their work.  The author thanks Bishop and Gauthier for sharing their work with him.

In the next section we fix some terminology and notation for clarity.  
In Section~\ref{the-proof} we present two lemmas and show that Theorem~\ref{theorem-generalized} follows from the lemmas.  The proofs of the lemmas are given in Sections~\ref{proof-of-lemma1} and~\ref{proof-of-lemma2}.

The arguments we give simplify somewhat when one aims to establish only Theorem~1.1 rather than the stronger Theorem~1.2.  Readers who are interested only in Theorem~1.1 may wish to consult the earlier version of this paper \cite{Izzo} available on the arXiv which deals exclusively with Theorem~1.1 and also discusses how, in that context at least, the use of dimension theory can be eliminated from the argument.

%
%
%
%

\section{Terminology and Notation}\label{notation}

We denote the Euclidean norm of a point $x$ in $\rn$ by $|x|$.
Given an $\rn$-valued map $f$ whose domain contains a subset $L$, we denote the restriction of $f$ to $L$ by $f|L$, and we set $ \|f\|_L = \sup\{ |f(x)| : x \in L \}$.
We denote the interior of a set $L$ by $L^\circ$.

We will use simplicial complexes in our proof.  
We will assume that the reader is familiar with the notion of a simplicial complex but recall the definition here to fix terminology and notation for clarity.  Readers for whom this material is unfamiliar may wish to consult for instance \cite{Munkres1966} and \cite{Munkres1984}.

Points $v_0,\ldots, v_n$ in $\R^N$ are said to be \emph{geometrically independent} if for real numbers $c_0,\ldots, c_n$ the equations $\sum_{j=0}^n c_j v_j=0$ and $\sum_{j=0}^n c_j=0$ imply $c_0 = \cdots = c_n =0$.
Given geometrically independent points $v_0,\ldots, v_n$ in $\R^N$, the \emph{$n$-simplex} $\s=v_0\cdots v_n$ they span is the set of points $x$ such that $x=\sum_{j=0}^n c_j v_j$, where $c_j\geq 0$ for all $j=0,\ldots, n$ and $\sum_{j=0}^n c_j=1$.  The points $v_0,\ldots, v_n$ are called the vertices of $\s$.
The number $n$ is called the \emph{dimension} of $\s$.
A \emph{face} of $\s$ is a simplex spanned by a subset of the vertices of $\s$.  The faces of $\s$ other than $\s$ itself are called the \emph{proper faces} of $\s$.  The \emph{boundary} $\partial \s$ of $\s$ is the union of the proper faces of $\s$.  The \emph{interior} $\s^\circ$ of $\s$ is defined by the equation $\s^\circ= \s \sm \partial \s$.

For our purposes, a \emph{simplicial complex} $K$ is a collection of simplices in $\R^N$ such that
\begin{enumerate}
\item[(1)] Every face of a simplex of $K$ is in $K$.
\item[(2)] The intersection of any two simplices of $K$ is a face of each of them.
\item[(3)] Each point belonging to a simplex of $K$ has a neighborhood in $\R^N$ that intersects only finitely many simplices of $K$.
\end{enumerate}
Condition (3) is often not included in the definition of a simplicial complex (and does \emph{not} follow from the other conditions), and furthermore, one can also consider simplicial complexes that are not contained in a finite-dimensional Euclidean space.  However, we will have no need of these more general objects, and so we adopt here the above more restrictive definition.  The \emph{dimension} of $K$ is defined to be the largest dimension of a simplex of $K$.  The union of the simplices of $K$ is called the \emph{polytope} of $K$ and denoted $|K|$.  We view $|K|$ as a topological space with the topology it inherits from $\R^N$.  It is easily verified that on account of condition (3) in our definition of simplicial complex, a subset $U$ of $|K|$ is open (closed) in $|K|$ if and only if $U\cap \s$ is open (closed) in $\s$ for every simplex $\s$ of $K$.  It follows that
a map on $|K|$ is continuous if and only if its restriction to each simplex of $|K|$ is continuous.

%
%
%
%

\section{Proof of Theorem~\ref{theorem-generalized}}\label{the-proof}

We will show that Theorem~\ref{theorem-generalized} follows from the following two lemmas.

\blem\label{lemma1}
Let $E$ be a subset of $\rn$, let $f:E\ra \rn$ be a continuous map that is zero-free on $E^\circ$, and let $\vep: E\ra\R_+$ be a positive continuous function.  Then there is a simplicial complex $K$ of dimension at most $n$, a continuous map $h: E\ra |K|$, a continuous map $g: |K| \ra \rn$ that is zero-free on each $n$-simplex of $K$, and a positive continuous function $\vep_K: |K|\ra \R_+$ such that for every $x\in E$ we have  $| f(x) -(g\circ h)(x) |<\vep(x)$ and $(\vep_K\circ h)(x)\leq \vep(x)$.
\elem

\blem\label{lemma2}
Let $X$ be a separable metrizable space, and let $g:X\ra\rn$ be continuous.  Suppose that $X$ is the union of two closed subspaces $A$ and $B$ such that $g|A$ is zero-free and $\dim B\leq n-1$.  Then $g$ can be approximated in the fine topology by a continuous map that is zero-free on all of $X$.
\elem

For the moment, we assume the lemmas and use them to prove Theorem~\ref{theorem-generalized}.  The proofs of Lemmas~\ref{lemma1} and~\ref{lemma2}  are given in Sections~\ref{proof-of-lemma1} and~\ref{proof-of-lemma2}, respectively.

\bpf[Proof of Theorem~\ref{theorem-generalized}]
\phantom{Step 1}

\noindent{\emph{Step 1: We prove the theorem in the special case $M=\rn$.}}

Fix a positive continuous function $\vep:E\ra\R_+$. 
By Lemma~\ref{lemma1}, there is a simplicial complex $K$ of dimension at most $n$, a continuous map $h: E\ra |K|$, a continuous map $g: |K| \ra \rn$ that is zero-free on each $n$-simplex of $K$, and a positive continuous function $\vep_K: |K|\ra \R_+$ such that $| f(x) -(g\circ h)(x)|<\vep(x)/2$ and $(\vep_K\circ h)(x)\leq \vep(x)$ for every $x\in E$.  Applying Lemma~\ref{lemma2}, with $X=|K|$, with $A$ the union of the $n$-simplices in $K$, and with $B$ the union of the simplices in $K$ of dimension at most $n-1$, yields a continuous map $\tg: |K| \ra \rn$ that is zero-free on all of $|K|$ such that $| g(x) - \tg(x) |<\vep_K(x)/2$ for every $x\in |K|$.  Then the continuous map $\tg\circ h : E \ra \rn$ is zero-free on $E$ and
for every $x\in E$ we have
\begin{align*}
| f(x) - (\tg \circ h)(x) | &\leq | f(x) - (g\circ h)(x)| + | (g\circ h)(x) - (\tg \circ h)(x) | \\
&< \bigl[\vep(x)/2\bigr] + \bigl[\vep_K(h(x))/2\bigr] \\
&\leq\vep(x).
\end{align*}
The special case of the theorem in which $M=\rn$ is proved.

\medskip
\noindent{\emph{Step 2: We prove the general case of the theorem.}}

Again fix a positive continuous function $\vep:E\ra\R_+$.  We may assume without loss of generality that $\vep$ is bounded above.

Since every component of a manifold is open, we may assume without loss of generality that $M$ is connected.  Then $M$ is second countable, and it follows that there exists a locally finite, countable collection of open sets $\{U_k\}_{k=1}^\infty$ of $M$ that covers $M$ and is such that for each $U_k$ there is a bounded open set $W_k$ of $\rn$ and a homeomorphism $\psi_k:U_k\ra W_k$.

Set $g_0=f$.
We will inductively construct a sequence $(g_j)_{j=1}^\infty$ of continuous maps from $E$ to $\rn$ such that for each $j=1, 2, \ldots$, the following conditions hold.
\begin{enumerate}
\item[(i)] $g_j=g_{j-1}$ on $E\sm U_j$.
\item[(ii)] $g_j$ is zero-free on $E\cap U_j$.
\item[(iii)] $g_j$ is zero-free on $E^\circ$.
\item[(iv)] $|g_j(x)-g_{j-1}(x)|< 2^{-j} \vep(x)$ for every $x\in E$.
\end{enumerate}

Before constructing the sequence $(g_j)_{j=1}^\infty$, we show that its existence will yield the theorem.  Condition~(iv), together with the boundedness of the function $\vep$,
insures that the sequence converges uniformly to a continuous map $g:E\ra \rn$, and furthermore, that for every $x\in E$ we have
$$
| f(x) - g(x) | = | g_0(x) - g(x) | \leq \sum_{j=1}^\infty |g_{j-1}(x) - g_j(x) | < \vep(x).
$$
Since the collection $\{U_k\}$ is a locally finite cover of $M$, given $x\in E$, there is a largest positive integer $m$ such that $x$ is in $U_m$.  Then $g_j(x)=g_m(x)$ for every $j\geq m$ by condition~(i), and hence $g(x)=g_m(x)$.  Since $g_m$ is zero-free on $E\cap U_m$ by condition~(ii), this gives $g(x)\not=0$.  Therefore, $g$ is zero-free on $E$.  Thus the proof will be complete once we establish the existence of the sequence $(g_j)_{j=1}^\infty$.

Let $k$ be a positive integer, and assume for the purpose of induction that we have chosen $g_j$ for all $1\leq j<k$.   (Note that when $k=1$, this assumption is vacuous.)  Set 
$\tE_k=\psi_k(E\cap U_k)\subset \rn$.  Set $\tg_k=g_{k-1}\circ(\psi_k^{-1}|\tE_k)$.  The interior of $\tE_k$ in $\rn$ is mapped by $\psi_k^{-1}$ into the interior of $E$, and hence $\tg_k$ is zero-free on the interior of $\tE_k$.  Denote the distance from a point $x\in W_k$ to the boundary $\partial W_k$ of $W_k$ in $\rn$ by $d(x, \partial W_k)$, and define $\tvep_k:\tE_k\ra \R_+$ by $\tvep(x) = \min\{ 2^{-k}(\vep\circ\psi_k^{-1})(x), d(x, \partial W_k)\}$.
By the special case of the theorem proved in Step~1, there exists a continuous map $\talpha_k:\tE_k\ra\rn$ such that 
$|\talpha_k(x)|< \tvep_k(x) {\rm\ for\ every\ } x\in \tE_k$ and $\tg_k + \talpha_k$ is zero-free on $\tE_k$.  Define a map $\alpha_k:E\ra \rn$ by
$$
\alpha_k(x) =
\begin{cases} 
\talpha_k(\psi_k(x)) &\mbox{for\ } x\in E\cap U_k \\
0& \mbox{for\ } x \in E\sm U_k.
\end{cases}
$$
Then $\alpha_k$ is continuous at every point of $E\cap U_k$ by the continuity of $\psi_k$ and $\talpha_k$.  Because $\tvep_k$ goes to zero at the boundary of $W_k$, for any constant $\eta>0$, there exists a compact set $C$ in $U_k$ such that $|\alpha_k(x)|<\eta$ for every $x\in E\sm C$, and hence $\alpha_k$ is continuous at every point of $E\sm U_k$.  Thus $\alpha_k$ is continuous.

Now set $g_k=g_{k-1} + \alpha_k$.  The inductive proof will be complete once we verify that with this definition of $g_k$, conditions~(i)--(iv) hold for $j=k$.  That $g_k=g_{k-1}$ on $E\sm U_k$ is obvious.  On $E\cap U_k$ we have
\begin{align*}
g_k|(E\cap U_k)  &= g_k\circ (\psi_k^{-1}|\tE_k) \circ (\psi_k|(E\cap U_k)) \\
&=\Bigl[g_{k-1} \circ (\psi_k^{-1}|\tE_k) + \alpha_k\circ(\psi_k^{-1}|\tE_k)\Bigr] \circ \Bigl[\psi_k|(E\cap U_k)\Bigr]\\
&=(\tg_k+\talpha_k) \circ (\psi_k|(E\cap U_k)).
\end{align*}
Thus $g_k$ is zero-free on $E\cap U_k$ because $\tg_k +\talpha_k$ is zero-free on $\tE_k$.  Since $g_{k-1}$ is zero-free on $E^\circ$, the conditions on $g_k$ just verified give that $g_k$ is also zero-free on $E^\circ$.  Finally, for every $x\in E\cap U_k$ we have
$$
|g_k(x) - g_{k-1}(x)|= |\alpha_k(x)|= |\talpha_k(\psi_k(x))|<|\tvep_k(\psi_k(x))|\leq 2^{-k}\vep(x),
$$
while (as already noted) for $x\in E\sm U_k$ we have $|g_k(x)-g_{k-1}(x)|=0$.
Thus conditions (i)--(iv) all hold with $j=k$.
\epf

%
%
%
%

\section{Proof of Lemma~\ref{lemma1}}\label{proof-of-lemma1}

The proof of Lemma~\ref{lemma1} will use the following lemma.

\blem\label{lemma3}
Let $E$ be a subset of $\rn$, let $\u$ be a collection of open sets of $\rn$ that covers $E$, and let $\Omega$ be the union of the members of $\u$.  There is an $n$-dimensional simplicial complex $L$ and a homeomorphism $\t : |L| \ra \Omega$ such that each simplex of $L$ is mapped by $\t$ into some member of $\u$ and for each $n$-simplex $\s$ in $L$ either $\t(\s) \subset E^\circ$ or else $\t(\s^\circ)$ intersects $\Omega\sm E$.
\elem

\bpf
By standard results about the existence of triangulations, there exists an $n$-dimensional simplicial complex $L$ and a homeomorphism $\tau_0:|L|\ra\Omega$ such that $\tau_0$ maps each simplex of $L$ into some member of $\u$.  (See for instance 
\cite[Theorem~10.6]{Munkres1966} and \cite[Theorem~16.4]{Munkres1984}.)  We will construct a homeomorphism $\gamma:\Omega\ra\Omega$ such that setting $\tau=\gamma\circ\tau_0$ yields a homeomorphism with all the properties asserted in the statement of the lemma.

Given a homeomorphism $\gamma:\Omega\ra\Omega$ and an $n$-simplex $\s\in L$, call $\s$ \emph{$\gamma$-good} if either $(\gamma\circ\tau_0)(\s) \subset E^\circ$ or $(\gamma\circ\tau_0)(\s^\circ)$ intersects $\Omega\sm E$, and call $\s$ \emph{$\gamma$-bad} otherwise.  Set $\gamma_0$ equal to the identity on $\Omega$ for notational convenience.  On account of condition (3) in our definition of simplicial complex, there are at most countably many simplices in $L$.    Thus there are at most countably many $\gamma_0$-bad $n$-simplices $\s_1, \s_2, \ldots$.  We will show that we can inductively choose homeomorphisms $\gamma_1, \gamma_2, \ldots$ such that for each $j=1, 2, \ldots$, every $n$-simplex that is $(\gamma_{j-1} \circ\cdots\circ\gamma_0)$-good is also $(\gamma_j \circ\cdots\circ\gamma_0)$-good, and in addition, $\s_j$ is 
$(\gamma_j \circ\cdots\circ\gamma_0)$-good, and $(\gamma_j\circ\cdots\circ\gamma_0)\circ \tau_0$ maps every simplex of $L$ into some member of $\u$.  In case there are only finitely many $\gamma_0$-bad $n$-simplices $\s_1,\ldots, \s_k$, then every $n$-simplex is $(\gamma_k\circ\cdots\circ\gamma_0)$-good, and setting $\tau=(\gamma_k\circ\cdots\circ\gamma_0)\circ\tau_0$ yields the lemma.  Otherwise, we will show $\gamma_1,\gamma_2,\ldots$ can be chosen in such a way that the sequence of homeomorphisms $(\gamma_j\circ\cdots\circ\gamma_0)_{j=0}^\infty$ converges uniformly to a homeomorphism $\gamma_\infty$ such that every $n$-simplex of $L$ is $\gamma_\infty$-good and $\gamma_\infty\circ\tau_0$ maps every simplex of $L$ into a member of $\u$.  Then setting $\tau=\gamma_\infty\circ\tau_0$ yields the lemma.

Now suppose that $\gamma:\Omega\ra\Omega$ is a homeomorphism such that $\gamma\circ\tau_0$ maps each simplex of $L$ into a member of $\u$, and suppose also that the inverse of $\gamma$ is uniformly continuous.  Consider an $n$-simplex $\s$ of $L$ that is $\gamma$-bad.  Then $(\gamma\circ\tau_0)(\s^\circ)\subset E$ and $\partial \s$ contains a point that is mapped by $\gamma\circ \tau_0$ to a point $p$ of $\partial E$.  Given any closed ball $B$ centered at $p$ contained in $\Omega$, there is a homeomorphism $\tgamma:\Omega\ra\Omega$ that is the identity outside of the interior of $B$ and sends a point of $(\gamma\circ\tau_0)(\s^\circ)$  to a point of $\Omega\sm E$.  Then $\s$ is $(\tgamma\circ\gamma)$-good.  Note that $\tgamma$ and $\tgamma^{-1}$ are each uniformly continuous.  Let $r$ denote the radius of $B$.  Then $\|(\tgamma\circ\gamma)-\gamma\|_{\Omega}\leq 2r$.  Note that for each $n$-simplex $\rho$ that is $\gamma$-good, if $\|(\tgamma\circ\gamma) - \gamma\|_\Omega$ is sufficiently small, then $\rho$ is also $(\tgamma\circ\gamma)$-good.  
Note also that for each simplex $\rho$, if $\|(\tgamma\circ\gamma) - \gamma\|_\Omega$ is sufficiently small, then $\rho$ is mapped by $(\tgamma\circ\gamma)\circ\tau_0$ into a member of $\u$.  
By choosing the radius $r>0$ of $B$ sufficiently small, we can arrange for $B$ to be disjoint from $\gamma\bigl(\tau_0(\rho)\bigr)$ for every simplex $\rho$ that does not share a vertex with $\s$; then $\tgamma$ is the identity on each such set $\gamma\bigl(\tau_0(\rho)\bigr)$.  Since there are only finitely many simplices that share a vertex with $\s$, it follows that for $r$ sufficiently small, every $\gamma$-good $n$-simplex is also $(\tgamma\circ\gamma)$-good
and 
every simplex of $L$ gets mapped by $(\tgamma\circ\gamma)\circ\tau_0$ into a member of $\u$.  Furthermore, because $\gamma^{-1}$ is assumed to be uniformly continuous, we can arrange to have $\|(\gamma^{-1}\circ\tgamma^{-1}) - \gamma^{-1}\|_\Omega$ smaller than any prescribed positive quantity by choosing $r$ small enough.  Note also that $\gamma^{-1}\circ\tgamma^{-1} = (\tgamma\circ\gamma)^{-1}$ is uniformly continuous.  Therefore, we can indeed inductively choose homeomorphisms $\gamma_1, \gamma_2, \ldots$ as in the preceding paragraph and such that in addition
\[
\|(\gamma_j\circ\cdots\circ\gamma_0) - (\gamma_{j-1}\circ\cdots\circ\gamma_0)\|_\Omega < 2^{-j}
\]
and
\[
\|(\gamma_j\circ\cdots\circ\gamma_0)^{-1} - (\gamma_{j-1}\circ\cdots\circ\gamma_0)^{-1}\|_\Omega < 2^{-j}.
\]
Then the two sequences $(\gamma_j\circ\cdots\circ\gamma_0)_{j=0}^\infty$ and
$((\gamma_j\circ\cdots\circ\gamma_0)^{-1})_{j=0}^\infty$ each converge uniformly to continuous functions $\gamma_\infty$ and $\gamma_{-\infty}$ which are inverses of each other and hence are homeomorphisms.  (Of course if there are only finitely many $\gamma_0$-bad simplices, then the sequences terminate and we do not need to consider convergence.)

Finally we need to note that every $n$-simplex of $L$ is $\gamma_\infty$-good and that every simplex of $L$ is mapped by $\gamma_\infty\circ\tau_0$ into a member of $\u$.  For that simply note that because each $\gamma_j$ is the identity on $(\gamma_{j-1}\circ\cdots\circ\gamma_0)\bigl(\tau_0(\rho)\bigr)$ for each simplex $\rho$ of $L$ that does not share a vertex with $\s_j$, the sequence $((\gamma_j\circ\cdots\circ\gamma_0)\circ\tau_0)_{j=0}^\infty$ eventually stabilizes on each simplex of $L$.
\epf

\bpf[Proof of Lemma~\ref{lemma1}]
Since $\vep$ is continuous on $E\subset\rn$, each point of $E$ has a neighborhood $W$ in $\rn$ such that the restriction of $\vep$ to $E\cap W$ is bounded away from zero.  It follows readily from this and the continuity of $f$ that there is a collection $\u$ of open balls of $\rn$ with centers in $E$ that covers $E$ such that for each ball $B\in\u$, the set $f(E\cap 5B)$ has diameter less than $\inf\{\vep(x): x\in E\cap 5B\}>0$, where $5B$ denotes the open ball whose center coincides with that of $B$ but whose radius is 5 times that of $B$.  Let $\Omega$ be the union of the members of $\u$.  By Lemma~\ref{lemma3}, there is an $n$-dimensional simplicial complex $L$ and a homeomorphism $\t : |L| \ra \Omega$ such that each simplex of $L$ is mapped by $\t$ into some member of $\u$ and for each $n$-simplex $\s$ in $L$ either $\t(\s) \subset E^\circ$ or else $\t(\s^\circ)$ intersects $\Omega\sm E$.  Henceforth, we will identify $\Omega$ with $|L|$ via the homeomorphism $\t$ and thus allow ourselves to regard $E$ as a subset of $|L|$, $\u$ as a collection of open subsets of $|L|$, etc.

Let $K$ be the simplicial complex obtained from $L$ by omitting those $n$-simplices of $L$ whose interiors intersect $\Omega\sm E$; in other words, $K$ consists of the $n$-simplices of $L$ contained in $E^\circ$ and the $(n-1)$-dimensional simplices forming the boundaries of the other $n$-simplices, and (of course) all faces of these simplices.  (Necessarily every point of $|L|$ lies in an $n$-simplex.)  For each $n$-simplex $\s$ of $L$ not in $K$, choose a point $z_\s$ in $\s^\circ\sm E$, and let $r_\s: \s \sm \{z_\s\} \ra \partial \s$ be radial projection.  The map $\th : |L|\sm \{z_\s : \s \in L \sm K\} \ra |K|$ defined by
$$
\th(z) =
\begin{cases} z &\mbox{for\ } z\in |K| \\
r_\s(z) & \mbox{for\ } z \in\s\sm \{z_\s\} \hbox{\ with\ } \s \in L\sm K
\end{cases}
$$
is continuous.  Let $h$ be the restriction of $\th$ to $E$.  

Let $Q$ be the union of the $n$-simplices in $K$.  Then $Q$ is a closed subset of $|L|$ and is contained in $E^\circ$.  For notational convenience, set $W=E^\circ$.  Set $\vv=\{B\sm Q : B\in \u {\rm\ and\ } B\sm Q \neq \emptyset\}$.    Then the collection $\vv \cup \{W\}$ is an open cover of $|L|$.  Choose a partition of unity $\{\pv\}_{V\in\vv\cup\{W\}}$ subordinate to $\vv\cup \{W\}$.  
For each $V\in \vv$, choose a ball $B \in \u$ such that $V=B\sm Q$, and let $z_V$ denote the center of $B$.  Define $\tg: |L| \ra \rn$ by
\[
\tg(z) = \varphi_W(z) f(z) + \sum_{V\in\vv} \pv(z) f(z_V)
\]
where of course we interpret $\varphi_W(z) f(z)$ to be zero for $z$ outside the support of $\varphi_W$.  Note that then $\tg$ is a well-defined continuous map.
Note also that $\tg=f$ on $Q\subset E^\circ=W$.  Thus $\tg$ is zero free on every $n$-simplex of $|K|$.

Consider an arbitrary $n$-simplex $\s$ of $L$ with $E \cap\s\not=\emptyset$.  The simplex $\s$ lies in some ball $B_0$ belonging to $\u$.  There are only finitely many sets $V\in\vv$ such that $\pv$ fails to be identically zero on $\s$, say $V_1,\ldots, V_m$.  Each $V_j$, $j=1,\ldots, m$, satisfies $V_j=B_j\sm Q$ for some ball $B_j\in \u$ with center $z_{V_j}$.  Let $B^*$ denote a ball among $B_0,\ldots, B_m$ with the largest radius.  Then because each of the balls $\row Bm$ intersects $B_0$, the ball $5B^*$ contains each of $B_0,\ldots, B_m$.  
Consequently, for every point $a\in E \cap \s$ we have $|f(a) -f(z_{V_j}) |< \inf\{\vep(x): x\in E \cap 5B^*\} \leq \vep(a)$ for each $j=1,\ldots, m$.
Thus for every point $a\in E \cap \s$ and $b\in\s$ we have
\begin{align*}
\bigl| f(a) - \tg(b) \bigr| &= \biggl | \Bigl[ \varphi_W(b) f(a) + \sum_{V\in\vv}\pv(b)  f(a) \Bigr] \\ &\qquad\qquad\qquad - \Bigl[ \varphi_W(b) f(b) + \sum_{V\in\vv}\pv(b)  f(z_V) \Bigr] 
\biggr | \\
&\leq \varphi_W(b)\bigl|f(a)-f(b)\bigr| + \sum_{V\in\vv} \pv(b) \bigl| f(a) -f(z_V)\bigr| \\
&<\varphi_W(b)\, \vep(a) + \sum_{V\in\vv} \pv(b)\, \vep(a) \\
&=\vep(a).
\end{align*}
Since $h(z)$ lies in $\s$ for every $z\in \s\cap E$, it follows at once that for every $z\in \s\cap E$ we have
\[
|f(z) - \tg(h(z))| < \vep(z).
\]
Therefore, the restriction of $\tg$ to $|K|$ has the properties required of the function $g$ in the lemma.

For each simplex $\s$ we have $\inf\{\vep(x):x\in E\cap\s\}>0$ because $\s$ is contained in a some member of $\u$.
Since the collection of simplices in $L$ is locally finite, it follows by a simple partition of unity argument left to reader (or see \cite[Theorem~41.8]{Munkres2000}) that there is a positive continuous function $\vep_L: |L|\ra \R_+$ such that for each simplex $\s\in L$ and each point $z\in\s$ we have $\vep_L(z)\leq \inf\{\vep(x):x\in E \cap\s\}$.  Taking $\vep_K$ to be the restriction of $\vep_L$ to $|K|$ completes the proof.
\epf

%
%
%
%

\section{Proof of Lemma~\ref{lemma2}}\label{proof-of-lemma2}

We will use a result from dimension theory whose statement involves the notion of two sets being separated by another set.
If $A_1$, $A_2$, and $B$ are pairwise disjoint subsets of a space $X$, we say that $A_1$ and $A_2$ are \emph{separated in $X$ by $B$} if $X\sm B$ can be split into two disjoint sets, open in $X\sm B$ and containing $A_1$ and $A_2$ respectively, i.e., if there exist sets $A'_1$ and $A'_2$, both open in $X\sm B$ such that
$$\openup1.25\jot\displaylines{
X\sm B = A'_1 \cup A'_2\cr
A_1\subset A'_1, \qquad A_2\subset A'_2\cr
A'_1\cap A'_2 = \emptyset.\cr}
$$

\bthm[\cite{Hurewicz-Wallman1948}, III 5 C]\label{separation}
Let $X$ be a separable metrizable space of dimension $\leq n-1$, and for each $i=1,\ldots, n$, let $C_i$ and $C'_i$ be closed subsets of $X$ such that $C_i\cap C'_i=\emptyset$.  Then there exist $n$ closed sets $\row Bn$ such that each $B_i$ separates $C_i$ and $C'_i$, and $B_1\cap\cdots\cap B_n=\emptyset$.
\ethm

The proof of the following result is a minor modification of the proof of 
\cite[Theorem VI 1]{Hurewicz-Wallman1948}, the special case in which the function $\vep$ is required to be constant, but we include it for the reader's convenience.

\bthm\label{unstable}
Let $X$ be a separable metrizable space of dimension $\leq n-1$, and let $f:X\ra \rn$ be continuous.  Let $\vep:X\ra \R_+$ be a positive continuous function.  Then there exists a continuous map $g:X\ra \rn$ such that $|f(x)-g(x)|<\vep(x)$ for every $x\in X$ and $g$ is zero-free.
\ethm

\bpf
Let $\row fn$ be the component functions of $f$.  Set 
\begin{eqnarray*}
\cpi=\{x\in X: f_i(x)\geq \vep(x)/2\}\phantom{-}\phantom{.} \\
\cmi=\{x\in X: f_i(x)\leq -\vep(x)/2\}.
\end{eqnarray*}
For each $i=1,\ldots, n$, the sets $\cpi$ and $\cmi$ are closed and disjoint.  Hence, by Theorem~\ref{separation} there exist closed sets $\row Bn$ such that $B_i$ separates $\cpi$ and $\cmi$, i.e., such that $X\sm B_i=U^+_i \cup U^-_i$, where $U^+_i$ and $U^-_i$ are disjoint open sets of $X$ that contain $\cpi$ and $\cmi$ respectively, and such that $B_1\cap \cdots\cap B_n=\emptyset$.  Define functions $\row gn:X\ra \rn$ by
$$
g_i(x) =
\begin{cases}
\vphantom{\vrule height 12pt depth 16pt}\phantom{-}f_i(x) & \mbox{if\ } x\in \cpi \cup \cmi\\
\vphantom{\vrule height 18pt depth 20pt}\phantom{-}\displaystyle{\frac{\vep(x)}{2} \frac {d(x, B_i)}{d(x, \cpi) + d(x, B_i)}}&\mbox{if\ } x\in U^+_i\sm \cpi\\
\vphantom{\vrule height 20pt depth 18pt}\displaystyle{-\frac{\vep(x)}{2} \frac {d(x, B_i)}{d(x, \cmi) + d(x, B_i)}}&\mbox{if\ } x\in U^-_i\sm \cmi\\
\phantom{-}0 &\mbox{if\ } x\in B_i.
\end{cases}
$$
One readily verifies that each $g_i$ is continuous and that $|g_i(x)-f_i(x)|<\vep(x)$ for every $x\in X$.  Moreover, $g_i(x)$ is zero only when $x$ is in $B_i$.  Thus since $B_1\cap\cdots\cap B_n=\emptyset$, the mapping $g=(\row gn)$ is zero-free.
\epf

We will also need two elementary lemmas.

\blem\label{Tietze}
Let $X$ be a normal space, let $A$ be a closed subspace of $X$, and let $f:A\ra \rn$ be continuous.  Let $m:X\ra \R_+$ be a positive continuous function such that $|f(x)|\leq m(x)$ for every $x\in A$.  Then there exists a continuous extension $g:X\ra \rn$ of $f$ such that $|g(x)|\leq m(x)$ for every $x\in X$.
\elem

\bpf
By the Tietze extension theorem, there is a continuous extension $\tf:X\ra \rn$ of $f$.  Define $g:X\ra \rn$ by
$$
g(x) =
\begin{cases} 
\vphantom{\vrule depth 11pt}\tf(x) &\mbox{if\ } |\tf(x)|\leq m(x) \\
\displaystyle\frac{m(x)}{|\tf(x)|}\, \tf(x) & \mbox{if\ } |\tf(x)|\geq m(x).
\end{cases}
$$
\epf

\blem\label{bounded-below}
Let $X$ be a space, and let $f:X\ra \rn$ be continuous and zero-free.  Let $\vep:X\ra \R_+$ be a positive continuous function.  Then there exists a continuous map $g:X\ra \rn$ such that for every $x\in X$ we have $|g(x)|\geq \vep(x)$ and
$|f(x)-g(x)|<\vep(x)$.
\elem

\bpf
Define $g:X\ra \rn$ by
$$
g(x) =
\begin{cases} 
\vphantom{\vrule depth 11pt}f(x) &\mbox{if\ } |f(x)|\geq \vep(x) \\
\displaystyle\frac{\vep(x)}{|f(x)|}\, f(x) & \mbox{if\ } |f(x)|\leq \vep(x).
\end{cases}
$$
\epf

\bpf[Proof of Lemma~\ref{lemma2}]
Fix a positive continuous function $\vep:X\ra\R_+$.  We are to find a continuous map $h: X\ra\rn$ that is zero-free and such that $| h(x)-g(x) | <\vep(x)$ for every $x\in X$.

Since $g|A$ is zero-free, Lemma~\ref{bounded-below} yields a continuous map $\tg:A\ra \rn$ such that for every $x\in A$ we have $|\tg(x)|\geq \vep(x)/2$ and $|\tg(x)-g(x)|<\vep(x)/2$.
By Lemma~\ref{Tietze}, there exists a continuous map $\alpha:X\ra\rn$ that extends $\tg-g|A$ and satisfies $|\alpha(x)|\leq \vep(x)/2$ for every $x\in X$.  On $A$ we have  $g+\alpha = \tg$.

Since $\dim B \leq n-1$, there is, by Theorem~\ref{unstable}, a continuous map $b:B\ra\rn$ such that $|b(x)|<\vep(x)/3$ for every $x\in B$ and the map $(g+\alpha)|B + b$ is zero-free.  By applying Lemma~\ref{Tietze} again, extend $b$ to a continuous map $\beta: X\ra \rn$ such that $|\beta(x)|\leq\vep(x)/3$ for every $x\in X$.  Set $h=g+\alpha+\beta$.  Then $h$ has the required properties.
\epf

\end{document}